\documentclass[11pt,twoside]{article} 
\usepackage{latexsym} 
\input{amssym.def} 
\input{amssym} 
\newfont{\fra}{eufm10 scaled 1095} 
\newfont{\Bb}{msbm10 scaled 1095} 
\newfont{\Bbg}{msbm10 scaled 1680} 
\newcommand\CC{{\mbox{\Bb C}}} 
\newcommand\RR{{\mbox{\Bb R}}}

\newcommand\fg{{\frak{g}}}

\newcommand\fl{{\frak l}}

\newcommand\fk{{\frak k}} 
\newcommand\fp{{\frak p}}

\newcommand\fa{{\frak a}}

\newcommand{\fso}{\mathop{{\frak s \frak o}}}

\newcommand{\Iso}{\mathop{{\rm Iso}}}

\newcommand{\Id}{{{\rm id}}} 
\newcommand{\ad}{{{\rm ad}}} 
\newcommand{\tr}{\mathop{{\rm tr}}}

\newcommand{\proj}{{{\rm pr}}} 
\newcommand\ip{{\langle\cdot \,,\cdot \rangle}}

\newcommand\proof{{\sl Proof. }} 
\newcommand{\qed}{\hspace*{\fill}\hbox{$\Box$}\vspace{2ex}} 
\newtheorem{theo}{Theorem}
\newtheorem{pr}{Proposition}

\newtheorem{re}{Remark}

\newtheorem{lm}{Lemma}
\begin{document} 
\title{Semisimplicity of indefinite extrinsic symmetric spaces and mean curvature} 
\author{Ines Kath}
\maketitle 
\begin{abstract} \noindent Improving a result of Eschenburg and Kim we give a criterion for semisimplicity of pseudo-Riemannian extrinsic symmetric spaces in terms of the shape operator with respect to the mean curvature vector.
\end{abstract}
\section{Introduction}
A non-degenerate submanifold $M$ of a pseudo-Euclidean space $\RR^{p,q}$ is called extrinsic symmetric if for each $x\in M$ the reflection $s_x:\RR^{p,q}\rightarrow \RR^{p,q}$ at the affine normal space maps $M$ to $M$. Extrinsic symmetric spaces are exactly those complete submanifolds of $\RR^{p,q}$ that have a parallel second fundamental form \cite{S}. 
A large class of examples can be obtained in the following way. Let $\fg$ be a Lie algebra and assume that $\fg$ admits an  $\ad(\fg)$-invariant non-degenerate inner product and an orthogonal Cartan decomposition $\fg=\fk\oplus \fp$, i.e., $[\fk,\fk]\subset\fk$, $[\fk,\fp]\subset\fp$ and $[\fp,\fp]=\fk$. Define $G:=\langle \exp \ad(X)|_\fp\mid X\in\fk\rangle$. Suppose that $\xi\in \fp$ satisfies $\ad(\xi)^3=-\ad(\xi)$. Then the orbit $G\cdot\xi$ is an extrinsic symmetric space in $\fp\cong \RR^{p,q}$. These examples will be called extrinsic symmetric spaces of Ferus type. This notion was introduced in \cite{EK} for the following reason.
In the case where the ambient space is Euclidean Ferus \cite{F2,F3} proved that each extrinsic symmetric space decomposes into a flat factor and a compact extrinsic symmetric  space and that, moreover, each compact extrinsic symmetric space $M\subset\RR^n$ arises in the way described above from a compact Lie algebra $\fg$. This leads to a complete classification in the case of a Euclidean ambient space. 

In the pseudo-Riemannian situation we cannot expect to get such a classification for arbitrary $(p,q)$ since, in general, a pseudo-Riemannian extrinsic symmetric space is not of Ferus type and even those that are cannot be classified since in `most cases' the Lie algebra $\fg$ is not reductive. However, as a first step one can determine all extrinsic symmetric spaces that arise by the above described construction from a semisimple Lie algebra $\fg$. This was done by Naitoh who gave a complete classification of these spaces~\cite{Nai}.
 
In \cite{EK}, Theorem B, Eschenburg and Kim gave a sufficient condition for being of Ferus type in terms of the shape operator $A$ and the mean  curvature vector $h$ of $M$. They proved that a full and indecomposable extrinsic symmetric space $M\subset \RR^{p,q}$ is of Ferus type if 
\begin{equation} \label{E*}
{A_h}^2\not=0.
\end{equation}
More exactly, for any full extrinsic symmetric space $M\subset V=\RR^{p,q}$ they constructed a metric Lie algebra  $(\fg,\ip)$ with an orthogonal Cartan decomposition $\fg=\fk\oplus V$. Under the Condition~(\ref{E*}) they find an element $\xi \in V$ such that $M$ is the $G$-orbit through $\xi \in V$ as explained above.

In \cite{Kext1} we showed that extrinsic symmetric spaces correspond bijectively to so-called weak extrinsic symmetric triples, which consist of a Lie algebra $\fg$ together with an invariant inner product, an involution $\theta$ and a derivation $D$ on $\fg$ satisfying certain conditions. If the extrinsic symmetric space is full then $\fg$ and $\theta$ coincide with the Lie algebra and the Cartan decomposition constructed in \cite{EK}. If, in addition, $M\subset\RR^{p,q}$ is of Ferus type, then  $D=\ad(\xi)$. 

We will  show that Condition~(\ref{E*}) is very strong. In fact, we will prove:
\begin{theo}\label {T}
An indecomposable extrinsic symmetric space satisfies (\ref{E*}) if and only if $\fg$ is semisimple.  
\end{theo}
Any extrinsic symmetric space that arises from a semisimple Lie algebra $\fg$ is of Ferus type. Hence Theorem 1 implies that  the indecomposable extrinsic symmetric spaces satisfying (\ref{E*}) are exactly those listed in \cite{Nai}.
Note that we do not use the  assumption of fullness, which is somewhat inadequate in the indefinite case.


I would like to thank Martin Olbrich for several useful discussions.

\section{Characterisation of semisimplicity by mean curvature}
Let $M\subset \RR^{p,q}$ be an $n$-dimensional extrinsic symmetric space. We denote by $\alpha$ the second fundamental form of $M$ and by  $A$ its shape operator. Let $V$ be the (possibly degenerate) smallest affine subspace that contains $M$. Then $M\subset V$ is a full extrinsic symmetric space. Obviously, $\alpha$ coincides with the second fundamental form of $M$ as a submanifold of $V$, which we also denote by $\alpha$. Furthermore, the shape operator $\bar A$  of $M\subset V$ satisfies $A_\eta=\bar A_\eta$ for each $\eta\in T^\perp M\cap V$. 

The vector $$h=\frac1n \tr \alpha\in V$$ is called mean curvature vector. In other words, if $x\in M$ and $e_1,\dots,e_n$ is an orthonormal basis of $T_xM$, then 
$$
h(x)=\frac 1n \sum_{j=1}^n \kappa_j \alpha(e_j, e_j),
$$
where $\kappa_j=\langle e_j,e_j\rangle=\pm1$.

 In \cite{Kext1} it is shown that the infinitesimal object associated to an extrinsic symmetric space $M\subset \RR^{p,q}$  is a so-called weak extrinsic symmetric triple. 
An extrinsic symmetric triple $(\fg,\Phi,\ip)$ consists of a metric Lie algebra $(\fg,\ip)$ and a pair $\Phi=(D,\theta)$,  where $\theta$ is an isometric involution on $\fg$ and $D$ is an antisymmetric derivation on $\fg$ such that $D\theta=-\theta D$ and $D^3=-D$.  Besides $\theta$ there is a further involution on $\fg$, namely $\tau=\exp \pi D$. We denote by $\fg_+$ and $\fg_-$ the eigenspaces of $\theta$, by $\fg^+$ and $\fg^-$ the eigenspaces of $\tau$ and by $\fg^+_+=\fg_+\cap\fg^+$, etc. the intersections of these eigenspaces. We require that $\fg_+^+=[\fg_+^-,\fg_+^-]$. A weak extrinsic symmetric triple is defined in the same way but the  inner product is allowed to be degenerate on $\fg^+_-$.   If $(\fg,\Phi,\ip)$ is a weak extrinsic symmetric triple, then its metric radical $R:=\fg\cap\fg^\perp\subset\fg_-^+$ is contained in the centre of $\fg$ and $\fg/R$ is an extrinsic symmetric triple. A weak extrinsic symmetric triple is called full if $\fg^+_-=[\fg_+^-,\fg_-^-]$.

There is a one-to-one correspondence between full extrinsic symmetric spaces (in a possibly degenerate ambient space) and full weak extrinsic symmetric triples. Let us briefly recall this correspondence, for details see \cite{Kext1}. 
As vector spaces $V=\fg_-$. Moreover, if $\phi$ is the map from $\fg_+$ to the Lie algebra $\frak{iso}(\fg_-)=\fso(\fg_-)\ltimes\fg_-$ of $\Iso(\fg_-)$ defined by
$$\phi(X)=((\ad X)|_{\fg_-},-D(X)),\quad X\in\fg_+,$$
then  $M$ is the orbit of $x_0=0$ under the action of the group $$\langle\, \exp(\phi(X))\mid X\in\fg_+\,\rangle \subset \Iso(\fg_-).$$ 
In particular, $T^\perp_{x_0}M=\fg_-^+$, $T_{x_0}M=\fg_-^-$.  With this identification the second fundamental form and the shape operator are given by
\begin{equation}\label{EdT} 
\alpha(u,v)=[Du,v],\quad A_\eta u =-[Du,\eta]
\end{equation}
for all $u,v\in \fg_-^-$ and $\eta\in\fg_-^+$. The Riemannian curvature tensor and the normal curvature satisfy
$$
R^M(u,v)w=-[[u,v],w],\quad R^\perp(u,v)\eta=-[[u,v],\eta]
$$
for $u,v,w\in \fg^-_-$ and $\eta\in \fg_-^+$.

\begin{re}{\rm
\begin{enumerate} 
\item If $M\subset V=\RR^{p,q}$ is full and of Ferus type, then $\fg=\fg_+\oplus\fg_-$ coincides with $\fg=\fk\oplus V$ as a metric Lie algebra with orthogonal Cartan decomposition. Moreover $D$ equals $\ad(\xi)$. The extrinsic symmetric spaces in both models are essentially the same, they just differ by a translation by $\xi$.
\item  The extrinsic symmetric triple $\fg/R$ is associated with the projection of the immersion of $M$ into $V$ to an immersion of $M$ into $V/(V\cap V^\perp)$ described in \cite{EK}, Section 2. 
\item The derivation $D$ belonging to a weak extrinsic symmetric triple is inner, i.e., $D=\ad(\xi)$ for some $\xi\in\fg_+^-$, if and only if the associated extrinsic symmetric space is of Ferus type. In particular, if $\fg$ is semisimple, then the associated extrinsic symmetric space is of Ferus type.
\end{enumerate}
}
\end{re}
An extrinsic symmetric space  $M\subset \RR^{p,q}$ is called indecomposable if the embedding $M\subset \RR^{p,q}$ does not decompose as a non-trivial direct product of embeddings $M_i\subset \RR^{p_i,q_i}$, $i=1,2$, with $M=M_1\times M_2$, $p=p_1+p_2$, $q=q_1+q_2$. Here we want to assume that $M\subset \RR^{p,q}$ is indecomposable. Then $(\fg,\ip,\Phi)$ is indecomposable, i.e., it is not the direct sum of non-trivial weak extrinsic symmetric triples. We have two kinds of indecomposable weak extrinsic symmetric triples, namely those that are semisimple and those that do not have simple ideals. We will show that the semisimple ones satify $A_h^2\not =0$ and that for those without simple ideals  $A_h^2=0$ holds. This will prove the theorem.

 \begin{pr}\label{P1}
If $\fg$ is semisimple, then $h=\lambda \xi$, where $\lambda\not=0$ is a complex number if $\fg$ is a complex Lie algebra (considered as a real one) or a real number if $\fg$ is not complex.
In particular, $A_h=-\lambda\cdot \Id$. 
\end{pr}
\proof
If $(\fg,\Phi,\ip)$ is indecomposable and $\fg$ is semisimple, we have two possibilities. Either $\fg$ is simple or $\fg$ is the sum of two isomorphic simple Lie algebras. In the latter case the two simple summands are orthogonal to each other and $\theta$ is an isometry between them. Hence,  $\langle X,Y\rangle=B_\fg(\mu X,Y)$, where $B_\fg$ is  the Killing form of the (real) Lie algebra $\fg$ and $\mu\not=0$ is a complex number if $\fg$ is a complex Lie algebra or a real number if $\fg$ is not complex. Let $e_1,\dots,e_n$ be an orthonormal basis of $\fg_-^-$. For any $\eta\in\fg_-^+$ we have
\begin{eqnarray*}
n\langle  h,\eta\rangle&=& \sum_{j=1}^n  \langle \kappa_j \alpha(e_j, e_j),\eta\rangle \ =\   \sum_{j=1}^n \langle \kappa_j [De_j, e_j],\eta\rangle\\
&=&\sum_{j=1}^n \kappa_j\langle   [[\xi,e_j], e_j],\eta\rangle\ =\   \sum_{j=1}^n \kappa_j \langle[\eta,[\xi,e_j]], e_j\rangle\\
&=& \frac12\Big( \sum_{j=1}^n \kappa_j \langle[\eta,[\xi,e_j]], e_j\rangle\ +\ \sum_{j=1}^n \kappa_j \langle D[\eta,[\xi,e_j]], D e_j\rangle\Big)\\
&=& \frac12\Big( \sum_{j=1}^n \kappa_j \langle[\eta,[\xi,e_j]], e_j\rangle\ +\ \sum_{j=1}^n \kappa_j \langle [\eta,[\xi,De_j]], D e_j\rangle\Big).
\end{eqnarray*}
The latter term equals $(1/2)\cdot \tr (\ad(\eta)\circ \ad(\xi))|_{\fg^-}$, which is equal to $(1/2)\cdot B_{\fg}(\eta,\xi)$ because of $\ad (\xi)|_{\fg^+}=0$. Hence $\langle  h,\eta\rangle= \langle\lambda \eta,\xi\rangle$ for $\lambda= 1/(2n\mu)$ and, consequently, $h=\lambda\xi$. The second assertion now follows from~(\ref{EdT}) and $\ad(\xi)^2=-\Id$ on $\fg^-$.\qed
\begin{re}{\rm
More exactly, if the Lie algebra $\fg$ is complex, then the involution $\theta$ can be complex linear or antilinear. In the first case, in general, $h$ is a complex multiple of $\xi$. In the second case $h$ must be a real multiple of $\xi$. }
\end{re}

\begin{pr} \label{P2}
If ${A_h}^2\not=0$, then $\fg$ is semisimple.
\end{pr}
We will give two proofs of this proposition. The first one combines results from \cite{EK} with the description of extrinsic symmetric spaces by weak symmetric triples. The second one is a short proof, which relies on the method of quadratic extensions.

{\sl First proof or Proposition \ref{P2}. } We will need the following Lemma. It is essentially due to Kim and Eschenburg, we just generalise it to the case of weak extrinsic symmetric triples.
\begin{lm}\label{L1} For any $u,v\in\fg_-^-$ we have 
$$B_\fg(u,v)=-2n \langle A_hu,v\rangle=B_\fg(Du,Dv).$$ 
\end{lm}
\proof   Recall that the metric radical $R:=\fg\cap\fg^\perp\subset \fg^+_-$  belongs to the centre of $\fg$ (cf.~\cite{Kext1}) and that the quotient $\bar\fg=\fg/R$ is an extrinsic symmetric triple. In particular, $\bar \fg=\fg_+^-\oplus \fg_+^-\oplus \fg_-^+/R\oplus\fg_-^-$ as a vector space. Let $\proj$ denote the projection from $\fg$ to $\bar \fg$. Since $R$ is in the centre of $\fg$ we have $B_\fg(u,v)= \tr (\ad(u)\circ \ad(v))=B_{\bar \fg}(u,v)$. Now we can apply \cite{EK}, Lemma 5.2., which says that the claim is true if $\fg$ is non-degenerate (for $u\in\fg_-^-$  the vector $t_u\in \fg_+^-$ used in \cite{EK} coincides with the vector  $Du$ in our notation). Consequently, $B_{\bar\fg}(u,v)=-2n\langle  A_{\bar h}u,v\rangle $ for 
\begin{equation}\label{Ebarh}\bar h:=(1/n)\cdot\sum _{j=1}^n \kappa_j [De_j,e_j]_{\bar\fg}=(1/n)\cdot\sum_{j=1}^n \kappa_j\proj\, [De_j,e_j]_{\fg}= \proj \,h,
\end{equation}
where $e_1,\dots, e_n$ is an orthonormal basis of $\fg_-^-$.
This implies
$$B_\fg(u,v)=B_{\bar\fg}(u,v)=-2n\langle  [\proj\, h,Du]_{\bar \fg},v\rangle_{\bar\fg}=-2n\langle  [h,Du],v\rangle = -2n \langle A_hu,v\rangle.$$
This also implies the second equation since every derivation is antisymmetric with respect to the Killing form and $D^2=-\Id$ holds on $\fg^-$.\qed

Also the next Lemma relies on ideas developed in \cite{EK}.
\begin{lm} \label{L2} Either $A_h$ is invertible or $A_h^2=0$.
\end{lm}
\proof 
Since $M\subset V$ is indecomposable we conclude from Moore's Lemma that $A_h$ has just one conjugate pair $\lambda,\bar\lambda$ of eigenvalues. If $\lambda\not=0$, then $A_h$ is invertible. Now suppose $\lambda=0$.

Since $h$ is parallel  $0=R^\perp(u,v)h= -[[u,v],h]$ for all $u,v\in\fg_-^-$. Since $\fg_+^+=[\fg_+^-,\fg_+^-]=[\fg_-^-,\fg_-^-]$ we obtain $[\fg_+^+,h]=0$, which also implies $[\fg_-^+,h]=0$ by $[\fg_-^+,h]\subset\fg_+^+$ and $\langle [\fg_-^+,h],\fg_+^+\rangle=\langle \fg_-^+,[\fg_+^+,h]\rangle=0$.
Take $u,v\in\fg_-^-$ and put $\eta:=[Du,v]$. Let $e_1,\dots,e_n$ be an orthonormal basis of $\fg_-^-$. Since $[h,\fg^+]=0$ we get
\begin{eqnarray*}
B_\fg(\eta,h) &=&\sum_{j=1}^n \kappa_j\Big( \langle \ad (\eta)\circ\ad (h) (e_j), e_j\rangle +\langle \ad (\eta)\circ\ad (h) (De_j), De_j\rangle\Big)\\
&=&2\sum_{j=1}^n \kappa_j \langle[\eta,[h,De_j]], De_j\rangle\ =\ 2\sum_{j=1}^n \kappa_j \langle[h,De_j], [\eta,De_j]\rangle\\
&=& 2\sum_{j=1}^n \kappa_j \langle A_he_j,A_\eta e_j\rangle =2\tr A_hA_\eta.
\end{eqnarray*}
Since $h$ is parallel the Ricci equation gives $0=\langle R^\perp(u,v)h,\eta\rangle=-\langle[A_h,A_\eta]u,v\rangle$ for all $u,v\in\fg_-^-$. Hence $[A_h,A_\eta]=0$, which implies that $A_hA_\eta$ is nilpotent, thus $\tr A_hA_\eta=0$.
Consequently, $B_\fg(\eta,h) =0$. On the other hand,
$$B_\fg(\eta,h)=B_\fg([Du,v],h)=-B_\fg(v,[Du,h])=B_\fg(v,A_hu)=c\langle A_hv,A_hu\rangle =c\langle A_h^2u,v\rangle$$
for some constant $c\not=0$ by Lemma~\ref{L1}. Hence $A_h^2=0$.\qed

By Lemma~\ref{L2} $A_h$ is invertible. Now Lemma~\ref{L1} implies that the Killing form of $\fg$ is non-degenerate on $\fg^-$.  Moreover, $[\fg^-,\fg^-]\supset [\fg^-_+,\fg_+^-]\oplus [\fg_+^-,\fg_-^-]=\fg_+^+\oplus\fg_-^+=\fg^+$, thus $[\fg^-,\fg^-]=\fg^+$. This implies that  the Killing form of $\fg$ is non-degenerate. Hence $\fg$ is semisimple.  \qed

{ \sl Second proof of Proposition \ref{P2}. }
Let $\fg$ be not semisimple. We have to show that $A_h^2=0$. By~(\ref{EdT}) it suffices to show that $(\ad h)^2=0$. This is equivalent to $(\ad \bar h)^2=0$ in $\bar\fg=\fg/R$, where $\bar h$ is defined as in~(\ref{Ebarh}). If $\fg$ is not semisimple, then also $\bar \fg$ is not semisimple. Hence we may assume that $\ip$ is non-degenerate. Since $(\fg,\ip,\Phi)$ is indecomposable $\fg$ does not contain simple ideals. Hence $\fg$ has the structure of a quadratic extension. This is proven in \cite{Kext2}. We do not want to recall this in detail here. We just collect the facts that we will need in the following.
As a vector space the Lie algebra $\fg$ can be identified with $\fl^*\oplus\fa\oplus\fl$ for some vector spaces $\fl$ and $\fa$ such that the Lie bracket on $\fg$ has the properties
\begin{equation}\label{Elb}
[\fa,\fa]\subset\fl^*,\ [\fl^*,\fl^*\oplus\fa]=0,\ [\fl^*,\fl]\subset\fl^*.
\end{equation}
Moreover, there are derivations $D_\fa:\fa\rightarrow \fa$, $D_\fl:\fl\rightarrow \fl$ and involutions $\theta_\fa:\fa\rightarrow\fa$, $\theta_\fl:\fl\rightarrow\fl$ such that  $D=(-D_\fl)^*\oplus D_\fa\oplus D_\fl$ and $\theta =\theta_\fl^*\oplus \theta_\fa\oplus \theta_\fl$. We define $\fa_-^-$, $\fl_-^-$ and $(\fl^*)^-_-$ analogously to $\fg_-^-$.  Let $a_1,\dots, a_k$ be an orthonormal basis of $\fa^-_-$ and $\kappa_j:=\langle a_j,a_j\rangle$. Furthermore, let $L_1,\dots ,L_l$ be a basis of $\fl^-_-$ and denote by $Z_1,\dots,Z_l$ its dual basis of $(\fl^*)^-_-$. Then $Z_1,\dots,Z_l,L_1,\dots ,L_l,a_1,\dots, a_k$ is a basis of $\fg_-^-=(\fl^*)^-_-\oplus \fa^-_-\oplus\fl^-_-$ and
$$
nh\ =\ \sum_{i=1}^l 2\alpha(L_i,Z_i)+\sum_{j=1}^k \kappa_j\alpha(a_j,a_j)   
\ =\ \sum_{i=1}^l 2[D L_i,Z_i]+\sum_{j=1}^k \kappa_j[Da_j,a_j].   
$$
Now (\ref{Elb}) shows that $h$ belongs to $\fl^*$. Using (\ref{Elb}) once more, this yields $(\ad\,h)^2=0$. 
\qed

We have seen that, for an indecomposable extrinsic symmetric space $M\subset \RR^{p,q}$, the operator $A_h$ is of one of the following three types:
\begin{enumerate}
\item $A_h$ is invertible. This holds if and only if $\fg$ is semisimple. More exactly, the following holds. If $\fg$ and $\theta$ are complex, then $A_h=-\lambda \Id$, where $\lambda\in\CC$, $\lambda\not=0$. Note that in this case $M$ as well as the ambient space have complex structures, thus this condition makes sense.  Semisimple extrinsic symmetric triples with complex $\fg$ and complex $\theta$ can be obtained as complexifications of the extrinsic symmetric triples that are associated with compact extrinsic symmetric spaces in Euclidean spaces. If $\fg$ is not complex or $\theta$ is not complex, then $A_h=-\lambda \Id$, where $\lambda\in\RR$, $\lambda\not=0$. 
A complete list of all semisimple extrinsic symmetric triples can be found in \cite{Nai}.
\item $A_h$ is 2-step nilpotent. Examples for this case can be found in \cite{Kext1}, Section 7 and \cite{Kext2}. For instance, all extrinsic symmetric embeddings of Cahen-Wallach spaces have a 2-step nilpotent operator $A_h$. 
\item $A_h=0$. If $M\subset\RR^{p,q}$ is full, this is equivalent to $h=0$.  Examples for this situation can also be found in \cite{Kext2}. E.g., the two-dimensional Lorentzian manifold
$$\RR^{1,1}\longrightarrow V,\quad (r,s) \longmapsto (r,s^2,s),$$
where $V\in\{\RR^{1,2}, \RR^{2,1}\}$, $\ip_V=2dx_1dx_3\pm dx_2^2$, is full and $h$ vanishes. The same is true for the three-dimensional Lorentzian manifold 
$$\RR^{1,2}\longrightarrow V,\quad (r,s,t) \longmapsto (s,-rt+r^4/4,t,r,r^2),$$
where $V= \RR^{2,3}$ with $\ip_V=2dx_1dx_4+2dx_2dx_5+ dx_3^2$.

If the minimal subspace containing $M$ is degenerate, then it can happen that $A_h=0$ but $h\not=0$. The following example arises as an extension of the above discussed embedding $\RR^{1,1}\hookrightarrow V$, $V\in\{\RR^{1,2}, \RR^{2,1}\}$.  Consider $\RR^5$ with the metric $\ip=2dx_1dx_3\pm dx_2^2+2dx_4dx_5$, which is isometric to the standard space $\RR^{2,3}$ in the `$+$'-case and to $\RR^{3,2}$ in the `$-$'-case. Let us denote this pseudo-Euclidean space by $W$. Then an easy computation shows that 
$$\RR^{1,1}\longrightarrow W, \quad (r,s)\longmapsto (r,s^2,s,rs\pm s^4/2,0)$$
is extrinsic symmetric with $h=(0,0,0,1,0)\not=0$ but $A_h=0$.
\end{enumerate}


\begin{thebibliography}{MMM} 
\bibitem[F1]{F2} D.\,Ferus, {\sl Immersions with Parallel Second Fundamental Form.} Math.\,Z. {\bf 140} (1974), 87 - 92. 
\bibitem[F2]{F3} D.\,Ferus, {\sl Symmetric Submanifolds of Euclidean Space.} Math.\,Ann.\,{\bf 247} (1980), 81 - 93. 
\bibitem[K1]{Kext1}I.\,Kath, {\sl Indefinite extrinsic symmetric spaces I.} J. reine angew.\,M., to appear.
\bibitem[K2]{Kext2}I.\,Kath, {\sl Indefinite extrinsic symmetric spaces II.} J. reine angew.\,M., to appear.
\bibitem[Ki]{Kim} J.\,R.\,Kim, {\sl Indefinite Extrinsic Symmetric Spaces.} Dissertation Augsburg (2005).
\bibitem[KE]{EK}  J.\,R.\,Kim, J.-H.\,Eschenburg,  {\sl Indefinite extrinsic symmetric spaces.} manuscr. math. {\bf 135} (2011), 203 - 214.
\bibitem[Nai]{Nai} H.\,Naitoh, {\sl Pseudo-Riemannian symmetric $R$-spaces.} Osaka J.\,Math. {\bf 21} (1984), 733 - 764.
\bibitem[S]{S} W.\,Str\"ubing, {\sl Symmetric submanifolds of Riemannian 
manifolds.} Math.\,Ann. {\bf 245} (1979), 37 - 44. 
\end{thebibliography}
\end{document}